\documentclass[a4paper]{article}
\usepackage[utf8]{inputenc}
\usepackage{amssymb,amsmath,latexsym,mathtools}
\usepackage[english]{babel}

\usepackage{tikz-cd}

\title{Field Generated by Division Points of Certain Formal Group Laws - II}
\author{Soumyadip Sahu}
\date{}

\begin{document}
\maketitle
\begin{abstract}
In this article we investigate several questions about the field generated by certain group laws following \cite{some results}.
\end{abstract}

\section{Introduction}
In \cite{some results} author studied some questions about the Galois group of field generated by division points of certain formal group laws and relation between this Galois group and the ring of endomorphisms of  corresponding formal group law. Present article is a continuation of this study and it's purpose is to fill up few obvious gaps of \cite{some results}.\\~\\
First, recall the set-up and main results of \cite{some results} : \\~\\
Let $p$ be a prime and let $K$ be a finite extension of $\mathbb{Q}_p$. Put $O_K$ to be the ring of integers of $K$, let $\mathfrak{p}_K$  
denote the unique maximal ideal of $O_K$ and $v_{\mathfrak{p}_K}(\cdot)$ be the valuation associated to it. Fix an algebraic closure $\overline{\mathbb{Q}}_p$ and $|.|_p$ be an fixed extension of the absolute value. Let $\overline{O}$ be the ring of integers of $\overline{\mathbb{Q}}_p$ and $\overline{\mathfrak{p}}$ be the unique maximal ideal of $\overline{O}$. Clearly $\overline{\mathfrak{p}} \cap K = \mathfrak{p}_K$. For $n \in \mathbb{N}$ let $\mu_{n}$ denote the group of $n$-th roots of unity inside $\overline{\mathbb{Q}}_p$.\footnote{In our convention $0 \notin \mathbb{N}$.}\\~\\
Let $A$ be an integrally closed, complete subring of $O_K$. Assume that $\mathfrak{p}_A$ is the maximal ideal, $K(A)$ is the field of fractions and $k(A)$ is the field of residues. Put $[k(A) : \mathbb{F}_p] = f_A$. Let $\mathfrak{F}$ be a (one dimensional, commutative) formal group-law defined over $O_K$ admitting an $A$ module structure. $\pi$ be a generator of $\mathfrak{p}_A$ and say \[[\pi](X) = \pi X + a_2X^2 + a_3X^3 + \cdots \in O_K[[X]] \tag{1.1}\] with at least one $a_i \in O_K - \mathfrak{p}_K$. Then $\min \,\{ i \,|\, |a_i|_p = 1 \} = p^h$ for some positive integer $h$ (see \cite[18.3.1]{haz2}). Now if $\pi_1$ is another generator of $\mathfrak{p}_A$ and \[[\pi_1](X) = \pi_1 X + b_2X^2 + \cdots \in O_K[[X]]\] then $b_{p^h} \in O_K - \mathfrak{p}_K$ and $\min \,\{ i \,|\, |b_i|_p = 1\} = p^h$. This integer $h$ is called the height of $\mathfrak{F}$ as formal $A$ module. If $a_i \in \mathfrak{p}_K$ for all $i \geq 2$ then we say, height of $\mathfrak{F}$ is infinity. We shall only consider formal $A$ modules of finite height.\\~\\
Let $A$ be as above and $\mathfrak{F}$ be a formal $A$ module over $O_K$. It defines a $A$ module structure on $\mathfrak{p}_K$ which naturally extends to a $A$-structure on $\overline{\mathfrak{p}}$. We shall denote the corresponding addition by $\oplus_{\mathfrak{F}}$ to distinguish it from usual addition.\\
Use $\text{End}_{O_K}(\mathfrak{F})$ to denote ring of endomorphisms of formal group $\mathfrak{F}$ which are defined over $O_K$ and $\text{End}_{O_K}^{A}(\mathfrak{F})$ to denote ring of endomorphisms of formal $A$ module $\mathfrak{F}$ which are defined over $O_K$. It is well-known that \[\text{End}_{O_K}(\mathfrak{F}) = \text{End}_{O_K}^{A}(\mathfrak{F})\;\; (\text{see}\,\cite[21.1.4]{haz2}).\]\\
Let $\pi$ be a generator of $\mathfrak{p}_A$. For each $n \geq 1$, use $\mathfrak{F}[\pi^n]$ to denote the $\pi^n$-torsion submodule of $\overline{\mathfrak{p}}$. For any sub-field $L$ of $\overline{\mathbb{Q}}_p$, let $L_{\mathfrak{F}}(\pi^n)$ be the subfield of $\overline{\mathbb{Q}}_p$ generated by $\mathfrak{F}[\pi^n]$ over $L$ and put 
\[ \begin{split}
\Lambda_{\pi}(\mathfrak{F}) = \bigcup_{i \geq 1} \mathfrak{F}[\pi^i],\\
L_{\mathfrak{F}}(\pi^{\infty}) = \bigcup_{i \geq 1} L_{\mathfrak{F}}(\pi^i).
\end{split}\]We shall adopt the convention $\mathfrak{F}[\pi^0] = \{0\}$.\\~\\
Let $\pi$ be a generator  of $\mathfrak{p}_K$ and put $K_{\pi} = \mathbb{Q}_p(\pi)$. Use $A_{\pi}$ to denote the ring of integers of $K_{\pi}$. For simplicity we shall write $\mathfrak{p}_{A_\mathfrak{\pi}}$ as $\mathfrak{p}_{\pi}$ and $f_{A_\pi}$ as $f_{\pi}$. Note that $\pi$ is a generator of $\mathfrak{p}_\pi$. Let $\mathfrak{F}$ be a formal $A_{\pi}$ module of height $h$ defined over $O_K$. We know that $f_{\pi}  \,|\, h$ (see \cite[Appendix-A, Proposition-1.1]{the}). Put $h_{r, \pi} = \frac {h} {f_{\pi}}$. If $\pi$ is clear from context we shall abbreviate $h_{r,\pi}$ as $h_r$. \\~\\ 
Fix a generator $\pi$ of $\mathfrak{p}_K$. Let $K_{\pi}^{h_r}/K_{\pi}$ be the unique unramified extension of degree $h_r$. Use $A_{\pi}^{h_r}$ to denote the ring of integers of this field.\\~\\
With this set-up the main results of \cite{some results} can be summarized as follows :\\~\\
\textbf{Remark 1.1 :} A) Let $\mathfrak{F}$ be a formal $A_{\pi}$ module of height $h$ defined over $O_K$. Assume that $\mu_{p^h - 1} \subseteq K$. Then the following are equivalent :\\
(i) $\text{Gal}(K_\mathfrak{F}(\pi^{\infty})|K)$ is abelian.\\
(ii) $K_{\mathfrak{F}}(\pi^n) = K(z)$ for all $z \in \mathfrak{F}[\pi^n] - \mathfrak{F}[\pi^{n-1}]$ and all $n \in \mathbb{N}$.\\
(iii) $\mathfrak{F}$ has a formal $A_{\pi}^{h_r}$ module structure defined over $O_K$ extending $A_{\pi}$ module structure. \\
B) Let $\mathfrak{F}$ be a formal $A_{\pi}$ module of height $h$ defined over $O_K$. Then, $\text{End}_{O_K}(\mathfrak{F})$ is integrally closed in its fraction field.\\~\\
Further, we have following conventions and observations :\\~\\
\textbf{Remark 1.2 :} i) Let $\pi$ be a generator of $\mathfrak{p}_K$ and let $\mathfrak{F}$ be a formal $A_{\pi}$ module over $O_K$. Such a formal module will be called a \emph{$\pi$-unramified group law} over $O_K$.\\
ii) Let $A$ be a complete, integrally closed subring of $O_K$ containing $A_{\pi}$. If $\mathfrak{F}$ is a formal $A$ module of height $h$ defined over $O_K$ then $f_A | h$ (see \cite[Remark - 1.2.v]{some results}). Note that $K(A)/K_{\pi}$ is unramified and $\pi$ is a generator of $\mathfrak{p}_A$. Using this argument one can put constraint on size of $\text{End}_{O_K}(\mathfrak{F})$. More precisely, if $F$ is fraction field of $\text{End}_{O_K}(\mathfrak{F})$ then $F/K_{\pi}$ is an unramified extension whose degree divides $h_{r}$. If the degree is exactly $h_r$ we say the endomorphism ring has \emph{full height}. (See \cite[Remark - 3.1.3]{some results})\\
iii) Let $L$ be an unramified extension of $K$ (possibly infinite). Use $\widehat{L}$ denote completion of $L$, $O_{\widehat{L}}$ for the ring of integers of $\widehat{L}$ and $\mathfrak{p}_{\widehat{L}}$ be the maximal ideal of $O_{\widehat{L}}$. If $\pi$ is a generator of $\mathfrak{p}_K$, then $\pi$ also generates $\mathfrak{p}_{\widehat{L}}$. Let $\mathfrak{F}$ be a $A_{\pi}$ module of height $h$ defined over $O_{\widehat{L}}$. Then one can prove analogues of results in remark-1.1 (note that in this case we shall replace the fields $K_{\mathfrak{F}}(\pi^n)$ by $\widehat{L}_{\mathfrak{F}}(\pi^n)$ and $K$ by $\widehat{L}$). (See \cite[Remark - 2.4.B]{some results})\\

\section{Field generated by division points}
Let $K$, $O_K$, $\pi$  be as in introduction. We introduce the following notations:\\
$E$ = a finite unramified extension of $K$,\\
$L = E^{\text{ur}} = K^{\text{ur}}$,\\
$\widehat{L}$ = completion of $L$,\\
$O_k$ = ring of integers of $k$ for any discrete valuation field $k$ ,\\
$\mathcal{F}_{\pi}(O_E)$ = set of all $A_\pi$ modules of finite height defined over $O_E$,\\
$\mathcal{F}_{\pi}(O_E)(h)$ = set of all $A_{\pi}$ modules of height $h$ defined over $O_E$,\\
$\mathcal{F}_{\pi}(O_{\widehat{L}})$ = set of all $A_{\pi}$ modules of finite height defined over $O_{\widehat{L}}$,\\
$\mathcal{F}_{\pi}(O_{\widehat{L}})(h)$ = set of all $A_{\pi}$ modules of height $h$ defined over $O_{\widehat{L}}$,\\
$C$ = completion of $\overline{\mathbb{Q}}_p$,\\
$\mathcal{E}_{E}$ = set of all subextensions of $\overline{\mathbb{Q}}_p/E$,\\
$\mathcal{E}_{L}$ = set of all subextensions of $\overline{\mathbb{Q}}_p/L$,\\
$\mathcal{E}_{\widehat{L}}$ = set of all subextensions of $C/\widehat{L}$.\\~\\
Now we have cannonical maps   
\[\begin{split}
D_{\pi, E} : \mathcal{F}_{\pi}(O_E) \to \mathcal{E}_E , \, \mathfrak{F} \to E_{\mathfrak{F}}(\pi^{\infty}) \\
D_{\pi}^{\text{ur}} : \mathcal{F}_{\pi}(O_L) \to \mathcal{E}_L , \, \mathfrak{F} \to L_{\mathfrak{F}}(\pi^{\infty})\\
\widehat{D}_{\pi} : \mathcal{F}_{\pi}(O_{\widehat{L}}) \to \mathcal{E}_{\widehat{L}}, \, \mathfrak{F} \to \widehat{L}_{\mathfrak{F}}(\pi^{\infty}).
\end{split}\]
The natural question that one can ask with with this set up is, when two elements have same image under $D_{\pi, E}$ ($D^{\text{ur}}_{\pi}$, $\widehat{D}_{\pi}$).\\
If $E$ is clear from context we shall drop it from subscript.\\~\\
\textbf{Remark 2.1 :} i) Let $\mathfrak{F}, \mathfrak{G} \in \mathcal{F}_{\pi}(O_E)$. Then $D_{\pi}(\mathfrak{F}) \subseteq D_{\pi}(\mathfrak{G}) \implies D_{\pi}^{\text{ur}}(\mathfrak{F}) \subseteq D_{\pi}^{\text{ur}}(\mathfrak{G}) \implies \widehat{D}_{\pi}(\mathfrak{F}) \subseteq \widehat{D}_{\pi}(\mathfrak{G})$.\\
Conversely, let $\mathfrak{F}, \mathfrak{G} \in \mathcal{F}_{\pi}(O_E)$. Then $\widehat{D}_{\pi}(\mathfrak{F}) \subseteq \widehat{D}_{\pi}(\mathfrak{G}) \implies D^{ur}_{\pi}(\mathfrak{F}) \subseteq D^{ur}_{\pi}(\mathfrak{G})$. The last statement can be proved by noting that $L(\mathfrak{F}[\pi^n])$ (resp. $L(\mathfrak{G}[\pi^n])$) is algebraic closure of $E$ in $\widehat{L}(\mathfrak{F}[\pi^n])$ (resp $\widehat{L}(\mathfrak{G}[\pi^n])$) for each $n \in \mathbb{N}$. \\
ii) Let $\mathfrak{F}, \mathfrak{G} \in \mathcal{F}_{\pi}(O_E)$ be two  group laws which are isomorphic over $O_E$. Then $D_{\pi}(\mathfrak{F}) = D_{\pi}(\mathfrak{G})$.\\
iii) Let $\mathfrak{F}, \mathfrak{G} \in \mathcal{F}_{\pi}(O_E)$  be such that there is a non-zero homomorphism of formal groups $f : \mathfrak{F} \to \mathfrak{G}$ over $O_E$. Since $\mathfrak{F}$ has finite height we conclude that $f$ is onto with finite kernel (see \cite[1.2]{lubin2}). Further by \cite[21.1.4]{haz2} $f$ is automatically a homomorphism of $A_{\pi}$ modules. Let $\beta \in \Lambda_{\pi}(\mathfrak{G})$. Then there is an $\alpha \in \Lambda_{\pi}(\mathfrak{F})$ such that $f(\alpha) = \beta$. Hence $\Lambda_{\pi}(\mathfrak{G}) \subseteq E_{\mathfrak{F}}(\pi^{\infty})$
and $E_{\mathfrak{G}}(\pi^{\infty}) \subseteq E_{\mathfrak{F}}(\pi^{\infty})$ ie $D_{\pi}(\mathfrak{G}) \subseteq D_{\pi}(\mathfrak{F})$.\\
Further, in such situation there is a non-zero homomorphism of $A_{\pi}$ modules $g : \mathfrak{G} \to \mathfrak{F}$ defined over $O_E$. (See \cite[1.6]{lubin2}) Since $\mathfrak{G}$ also has finite height we conclude $D_{\pi}(\mathfrak{G}) \subseteq D_{\pi}(\mathfrak{F})$. \\
Thus $D_{\pi}(\mathfrak{F}) = D_{\pi}(\mathfrak{G})$.\footnote{Note that our definition of height differs from Lubin's definition of height, but notion of finiteness is same.}\\ 
iv) Note that the arguments of $(ii)$ and $(iii)$ also hold for $D^{\text{ur}}_{\pi}$, $\widehat{D}_\pi$.\\~\\
In literature one usually considers arbitrary formal group laws (possibly without any unramified group law structure) and work with the field generated by $p$-torsion points. In such situation height means height as formal $\mathbb{Z}_p$ module and there is a notion of an absolute endomorphism ring (see \cite{lubin1}). We shall modify notations introduced above and use $p$ instead of $\pi$ in subscript, to change the set-up to  arbitrary formal group laws and field generated by $p$-torsions. Note that any homomorphism of formal groups is compatible with module structure provided both groups has that module structure. We can compare this general framework with our restricted situation as follows :\\~\\
\textbf{Remark 2.2 :} i) Similar implications as in remark-2.1(i).\\
ii) Let $\mathfrak{F}, \mathfrak {G}\in \mathcal{F}_p(O_E)$ be such that $\text{Hom}_{O_E}(\mathfrak{F}, \mathfrak{G}) \neq 0$. Then by similar arguments as in remark-2.1(iii) $D_p(\mathfrak{F}) = D_p(\mathfrak{G})$.\\
iii) Let $\mathfrak{F} \in \mathcal{F}_p(O_E)$. Note that $\text{End}_{O_E}(\mathfrak{F})$ is a closed subring of $O_E$ via the $c$ map (see \cite[2.2.1]{lubin1}). There is a $\mathfrak{G} \in \mathcal{F}_p(O_E)$ such that $\text{Hom}_{O_E}(\mathfrak{F}, \mathfrak{G}) \neq 0$, fraction field of $\text{End}_{O_E}(\mathfrak{F})$ = fraction field of $\text{End}_{O_E}(\mathfrak{G})$ (say, $F$) and $\text{End}_{O_E}(\mathfrak{G})$ is integrally closed (see \cite[3.2]{lubin2}). If $E/F$ is unramified, then $\mathfrak{G}$ is an unramified group law and we can use our set-up in questions related to field generated by $p$-torsion points of $\mathfrak{F}$.\\
iv) $\mathfrak{F} \in \mathcal{F}_{\pi}(O_E)(h) \implies \mathfrak{F} \in \mathcal{F}_p(O_E)(eh)$ where $e$ is ramification index of $K/\mathbb{Q}_p$. Further, $D_{\pi}(\mathfrak{F}) = D_{p}(\mathfrak{F})$ (similar statements for $D^{\text{ur}}$, $\widehat{D}$). \\~\\
Let $\mathfrak{F} \in \mathcal{F}_p(O_E)$. The absolute endomorphism ring of $\mathfrak{F}$ (denoted $\text{End}(\mathfrak{F})$) is defined as \[ \text{End}(\mathfrak{F}) = \bigcup_{k} \text{End}_{O_k}(\mathfrak{F}) \]
 where $k$ any complete discrete valuation field containing $E$ and $O_k$ is the associated ring of integers.\\
Assume that height of $\mathfrak{F}$ as $\mathbb{Z}_p$ module is $h$ and use $\mathcal{S}$ to denote all the field extensions of $\mathbb{Q}_p$ of degree $h$. Note that $\mathcal{S}$ is a finite set. Put $k_1$ to be the compositum of all elements in $\mathcal{S}$ and $E$. Clearly, $k_1/E$ is a finite extension. One can show that $\text{End}(\mathfrak{F}) = \text{End}_{O_{k_1}}(\mathfrak{F})$. Further $c(\text{End}_{O_k}(\mathfrak{F})) = O_k \cap c(\text{End}(\mathfrak{F}))$ for any $k$ as above. (See \cite[2.3.3]{lubin1})\\
It's easy to see that $\text{End}(\mathfrak{F})$ is a noetherian, local ring and its fraction field is a finite extension of $\mathbb{Q}_p$. \\~\\
Now we have the following proposition :\\~\\
\textbf{Proposition 2.3 :} Let $\mathfrak{F} \in \mathcal{F}_{\pi}(O_E)$. Then $\text{End}(\mathfrak{F})$ is integrally closed in its fraction field.\\~\\
\textbf{Proof :} Similar to proof of theorem-1.3 (see \cite[Theorem-4.1.2]{some results}). Only thing to notice in this case is $\mathfrak{m}$, the unique maximal ideal of $\text{End}(\mathfrak{F})$ is stable under action of $G_{E}$, the absolute Galois group of $E$. $\square$\\~\\
Let $\mathfrak{F}, \mathfrak{G} \in \mathcal{F}_{p}(O_{\widehat{L}})$. Use $\mathfrak{F}_p, \mathfrak{G}_p$ to denote the associated $p$-divisible groups arising from $p$-torsion points (see \cite[2.2]{tate}). Assume that $T$ is a finite extension of $\widehat{L}$ and $O_T$ is the ring of integers. Clearly $T/\widehat{L}$ is totally ramified. Now by a result of Waterhouse \[ \text{Hom}(\mathfrak{F}_p, \mathfrak{G}_p) = \text{Hom}((\mathfrak{F}_p)_{O_T}, (\mathfrak{G}_p)_{O_T})) \;\; \cite[\text{Theorem} \,3.2]{water}\]   
From theory of $p$-divisible groups (\cite[2.2]{tate}) we know the canonical maps 
\[\begin{split}
\text{Hom}_{O_{\widehat{L}}}(\mathfrak{F}, \mathfrak{G}) \to \text{Hom}(\mathfrak{F}_p, \mathfrak{G}_p),\\
\text{Hom}_{O_T}(\mathfrak{F}, \mathfrak{G}) \to \text{Hom}((\mathfrak{F}_p)_{O_T}, (\mathfrak{G}_p)_{O_T})\end{split}\]  
are bijective.\\
Note that $\text{Hom}_{O_{\widehat{L}}}(\mathfrak{F}, \mathfrak{G}) \subseteq \text{Hom}_{O_T}(\mathfrak{F}, \mathfrak{G})$. In light of discussion above, we must have $\text{Hom}_{O_{\widehat{L}}}(\mathfrak{F}, \mathfrak{G}) = \text{Hom}_{O_T}(\mathfrak{F}, \mathfrak{G})$.\\~\\
\textbf{Lemma 2.4 :} Let $\mathfrak{F} \in \mathcal{F}_p(O_{\widehat{L}})$. Then $\text{End}(\mathfrak{F}) = \text{End}_{O_{\widehat{L}}}(\mathfrak{F})$.\\~\\
\textbf{Proof :} Let $h$ be height of $\mathfrak{F}$ as $\mathbb{Z}_p$ module and $\mathcal{S}$ be as before. Use $T$ to denote compositum of elements of $\mathcal{S}$ and $\widehat{L}$. Clearly $T/\widehat{L}$ is a finite totally ramified extension. Now the lemma follows from discussion above. $\square$\\~\\
In what follows we shall identify the endomorphism ring with its image under $c$ map.\\~\\
 \textbf{Corollary 2.5 :} Let $\mathfrak{F} \in \mathcal{F}_p(O_E)$. Use $F$ to denote the fraction field of $\text{End}(\mathfrak{F})$. Then, $EF/E$ is a finite unramified extension.\\~\\
\textbf{Proof :}  $\text{End}_{O_{\widehat{L}}}(\mathfrak{F}) = \text{End}(\mathfrak{F}) \cap O_{\widehat{L}}$. From lemma-2.4 we have $\text{End}(\mathfrak{F}) \subseteq O_{\widehat{L}}$. Since $F/\mathbb{Q}_p$ is finite, $EF/E$ is a finite unramified extension. $\square$\\~\\
\textbf{Lemma 2.6 :} Let $\mathfrak{F} \in \mathcal{F}_{\pi}(O_E)(h)$. Then $\text{End}(\mathfrak{F})$ is an integrally closed subring of $A_{\pi}^{h_r}$ containing $A_{\pi}$.\\~\\
\textbf{Proof :} $F$ be as before. $EF/E$ is unramified extension. Note that, $E/K_{\pi}$ is unramified. Hence $EF/K_{\pi}$ is unramified. \\
By hypothesis, $K_{\pi} \subseteq F$. Thus $F/K_{\pi}$ is unramified. By proposition-2.3 $\text{End}(\mathfrak{F}) = O_F$. Assume $[k(O_F) : \mathbb{F}_p ] = f_{F}$. Then $f_{F} | h$. Hence $K_{\pi} \subseteq F \subseteq K_{\pi}^{h_r}$. The lemma follows from here. $\square$ \\~\\
\textbf{Corollary 2.7 :} Let $\mathfrak{F} \in \mathcal{F}_{\pi}(O_E)(h)$ be such that $\text{End}_{O_E}(\mathfrak{F}) = A_{\pi}^{h_r}$ ie the endomorphism ring has full height. Then, $\text{End}_{O_E}(\mathfrak{F}) = \text{End}(\mathfrak{F})$. \\~\\
\textbf{Proof :} Follows from lemma-2.6. $\square$\\~\\
\textbf{Corollary 2.8 :} Let $\mathfrak{F} \in \mathcal{F}_{\pi}(O_E)(h)$. Assume that $\mu_{p^h - 1} \subseteq O_E$ ie $h | f_E$ where $f_E$ is the degree of residue extension corresponding to $O_E$. Then $\text{End}(\mathfrak{F}) = \text{End}_{O_E}(\mathfrak{F})$. \\~\\
\textbf{Proof :} We have $\text{End}_{O_E}(\mathfrak{F}) = \text{End}(\mathfrak{F}) \cap O_E$. Now the result follows from lemma-2.6 and hypothesis. $\square$ \\~\\
\textbf{Remark 2.9 :} 
If $\text{End}_{O_E}(\mathfrak{F}) = A_{\pi}^{h_r}$, then $A_{\pi}^{h_r} \subseteq O_E$. So corollary-2.7 is a special case of corollary-2.8. \\~\\
Let $\mathfrak{F} \in \mathfrak{F}_p(O_E)(h)$. In Lubin's terminology (\cite[4.3.1]{lubin1}) $\mathfrak{F}$ is said to be full if the following holds :\\
1. $\text{End}(\mathfrak{F})$ is integrally closed in its fraction field,\\
2. $\text{End}(\mathfrak{F})$ is a free $\mathbb{Z}_p$ module of rank $h$.\\~\\ 
\textbf{Remark 2.10 :} If $\mathfrak{F}$ is an unramified group law of height $h$ whose endomorphism ring is full, it is full in the sense mentioned above (recall, $A_{\pi}^{h_r}$ is a free $\mathbb{Z}_p$ module of rank $eh$ and  where $e$ is ramification index of $K/\mathbb{Q}_p$ use remark-2.2(iv)).\\~\\
We have the following theorem :\\~\\
\textbf{Theorem 2.11 :} Let $O$ be a complete discrete discrete valuation ring and assume that $\mathfrak{F}$ and $\mathfrak{G}$ are formal group laws over $O$, which are $\mathbb{Z}_p$ modules of finite height. Suppose that $O$ is large enough so that $\text{End}_{O}(\mathfrak{F}) = \text{End}(\mathfrak{F})$, $\text{End}_{O}(\mathfrak{G}) = \text{End}(\mathfrak{G})$ and its residue field is algebraically closed. Further assume $c(\text{End}(\mathfrak{F})) = c(\text{End}(\mathfrak{G}))$. If $\mathfrak{F}$, $\mathfrak{G}$ are full group laws (in sense defined above) then they are isomorphic over $O$. \\~\\
\textbf{Proof :} See \cite[4.3.2]{lubin1}. $\square$ \\~\\
\textbf{Corollary 2.12 :} Let $\mathfrak{F},\mathfrak{G} \in \mathcal{F}_{\pi}(O_E)(h)$ with $\text{End}_{O_E}(\mathfrak{F}) = \text{End}_{O_E}(\mathfrak{G}) = A_{\pi}^{h_r}$. Then $\mathfrak{F}$ and $\mathfrak{G}$ are isomorphic over $O_{\widehat{L}}$. \\~\\
\textbf{Proof :} Put $O = O_{\widehat{L}}$. Now the corollary follows from lemma-2.4 and theorem-2.11. $\square$\\~\\
\textbf{Corollary 2.13 :} Let $\mathfrak{F}, \mathfrak{G}$ be as in statement of corollary-2.12. Then $\widehat{D}_{\pi}(\mathfrak{F}) = \widehat{D}_{\pi}(\mathfrak{G})$.\\~\\
\textbf{Proof :} Follows from corollary-2.12. $\square$\\~\\
\textbf{Corollary 2.14 :} Let $\mathfrak{F} \in \mathcal{F}_{\pi}(O_E)(h)$ with $\text{End}_{O_E}(\mathfrak{F}) = A_{\pi}^{h_r}$. Then $D^{\text{ur}}_{\pi}(\mathfrak{F}) = E^{\text{ab}}$ where $E^{\text{ab}}$ is the maximal abelian extension of $E$.\\~\\
\textbf{Proof :} Let $\mathfrak{G}$ be a Lubin-Tate module with respect to the parameter $\pi$ on $A_{\pi}^{h_r}$. Note that, the hypothesis on $\mathfrak{F}$ implies $A_{\pi}^{h_r} \subseteq O_{E}$. Clearly $\mathfrak{G} \in \mathcal{F}_{\pi}(O_E)(h)$ and $\text{End}_{O_E}(\mathfrak{G}) = A_{\pi}^{h_r}$ (the last part follows from lemma-2.6 and definition of Lubin-Tate module). Since $E/K_{\pi}^{h_r}$ is unramified, by class field theory we conclude $D_{\pi}^{\text{ur}}(\mathfrak{G}) = E^{\text{ab}}$ (see \cite[Chapter-III, Corollary-7.7]{neu}).\\
Now the result follows from corollary-2.13 and remark-2.1(i). $\square$\\~\\
\textbf{Remark 2.15 :} From corollary-2.14 one easily sees that $D_{\pi}^{\text{ur}}(\mathfrak{F})$ is same for any unramified group law over $O_E$ whose endomorphism ring has full height ie it is independent of $\pi$, $h$ and $\mathfrak{F}$. This along with remark-2.1(i) answers two questions posed in \cite[Remark-4.1]{some results}.

\section{The associated Galois representation}
In this section we recall some standard results about Galois representation arising from formal group laws and use them in context of unramified group laws.\\
Let $E$ be as before. Put 
\[\begin{split}
G_{E} = \text{Gal}(\overline{\mathbb{Q}}_p|E),\\
G_{L} = \text{Gal}(\overline{\mathbb{Q}}_p|L),\\
G_{\widehat{L}} = \text{Gal}(C|\widehat{L}).
\end{split}\] 
Note that $G_L = G_{\widehat{L}}$ and $G_L$ is a closed normal subgroup of $G_E$.\\
Let $\mathfrak{F} \in \mathcal{F}_p(O_E)(h)$. It is well-known that the $p$-adic Tate module $T_{p}(\mathfrak{F})$ is a free $\mathbb{Z}_p$ module of rank $h$ (\cite[1.2]{lubin2}). Put $V_p(\mathfrak{F}) = T_p(\mathfrak{F})\underset{\mathbb{Z}_p}\otimes\mathbb{Q}_p.$ It is a $\mathbb{Q}_p$ vector space of dimension $h$. Fixing an ordered base for $T_p(\mathfrak{F})$ one obtains natural representations :
\[\begin{split}
\rho_{p}(\mathfrak{F}) : G_E \to \text{Gl}_{h}(\mathbb{Q}_p),\\
\widehat{\rho}_{p}(\mathfrak{F}) : G_{\widehat{L}} \to \text{Gl}_h(\mathbb{Q}_p).
\end{split}\]\\
\textbf{Remark 3.1 :} i) The maps $\rho_p(\mathfrak{F})$ and $\widehat{\rho}_p(\mathfrak{F})$ factor through canonical inclusion $\text{Gl}_h(\mathbb{Z}_p) \to \text{Gl}_h(\mathbb{Q}_p)$.\\
ii) These maps are continuous with respect to usual pro-finite topologies on both sides.\\
iii) Put $\rho_{p}(\mathfrak{F})(G_E) = H$ and $\widehat{\rho}_{p}(\mathfrak{F})(G_{\widehat{L}}) = \widehat{H}$. Note that $G_{\widehat{L}}$ is a closed sub-group of $G_E$ and $G_E$ is compact. Further, $\text{Gl}_{h}(\mathbb{Q}_p)$ is Hausdorff. Hence $H$ and $\widehat{H}$ are closed subgroups of $\text{Gl}_h(\mathbb{Q}_p)$. By a $p$-adic analogue of Cartan's theorem we know, any closed sub-group of $\text{Gl}_{h}(\mathbb{Q}_p)$ is a $p$-adic Lie group in analytic sense. (see \cite[V.9, Corollary of Threorem-1]{serre 1})\footnote{In this article $p$-adic Lie group means analytic group of finite dimension over a finite extension of $\mathbb{Q}_p$. If the base field is not explicitly mentioned it is assumed to be $\mathbb{Q}_p$. }. So $H$ and $\widehat{H}$ are $p$-adic Lie groups. \\
iv) $\rho_p(\mathfrak{F})$ induces an isomorphism of abstract groups between $\text{Gal}(D_{\pi}(\mathfrak{F})|E)$ and $H$. Since it is continuous, $\text{Gal}(D_{\pi}(\mathfrak{F})|E)$ is compact and $H$ is Hausdorff this isomorphism of abstract group is actually an isomorphism of topological groups.\\
Similarly $\widehat{\rho}_p(\mathfrak{F})$ induces an isomorphism of abstract groups between $\text{Gal}(\widehat{D}_{\pi}(\mathfrak{F})|\widehat{L})$ and $\widehat{H}$.\\
v) In the following discussion we shall be mostly concerned about $\rho_p(\mathfrak{F})$. Properties of $\widehat{\rho}_p(\mathfrak{F})$ are similar and we shall mention them only if we need them.\\~\\
Next part of discussion follow works of Serre and Sen. (\cite{serre 2}, \cite{sen}).\\
Put $V_C = V_p(\mathfrak{F})\underset{{\mathbb{Q}}_p}{\otimes}C$. It is a $C$ vector space of dimension $h$ on which $G_E$ acts semi-linearly ie \[ s(cx) = s(c)s(x) \] for all $s \in G_E, c \in C, x \in V_C$.\\
Let $U_p$ be the group of units of $\mathbb{Z}_p$. Fix a generator $t$ for $T_p(\mathbb{G}_m)$ and let $\chi : G_E \to U_p$ be the corresponding cyclotomic character. For $i \in \mathbb{Z}$ define, \[ V^{i} = \{ x \in V_C \,|\, s(x) = \chi(s)^i x \; \forall s \in G_E \}. \]
It is a $E$ vector-space. Put $V(i) = C$-space spanned by $V^{i}$. Following theorem is a fundamental result of the theory : \\~\\
\textbf{Theorem 3.2 (Hodge-Tate decomposition)} :
\[V_C = V(0) \oplus V(1)\]      
as $C$-vector spaces. Further, $\text{dim}_C(V(0)) = h - 1$ and $\text{dim}_C(V(1)) = 1$.\\~\\
\textbf{Proof :} See \cite[Section-4, corollary-2 of theorem-3]{tate} and \cite[Section-5, Proposition-6]{serre 2}. $\square$ \\~\\
Let $H$ be as before. Use $\mathfrak{h}$ to denote the Lie algebra associated with $H$. Remember that the exponential map defines an isomorphism from a neighbourhood of $0$ in $\mathfrak{h}$ onto an open subgroup of $H$.\\
Use $H_{\text{alg}}$ to denote the smallest algebraic subgroup of $\text{Gl}_h(\mathbb{Q}_p)$ containing $H$. Similarly $\mathfrak{h}_{\text{alg}}$ be the smallest algebraic Lie sub-algebra of $\mathfrak{gl}(h, \mathbb{Q}_p)$ containing $\mathfrak{h}$. One can show that $\mathfrak{h}_{\text{alg}}$ is indeed the Lie algebra corresponding to $H_{\text{alg}}$ (\cite[Section-1, Proposition-2]{serre 2}). \\
As a consequence of theorem-3.2 the Galois representation in concern is a `Hodge-Tate' representation. The following result is due to Serre and Sen :\\~\\
\textbf{Proposition 3.3 :} Notations be as above. Then, \[\mathfrak{h}_{\text{alg}} = \mathfrak{h}\] ie $\mathfrak{h}$ is an algebraic Lie algebra. \\~\\
\textbf{Proof :} See \cite[Section-6, Theorem-2]{sen}. $\square$\\~\\
 \textbf{Remark 3.4 :} Proposition-3.3 answers lie algebra version of a question in \cite[Remark-3.2.1]{some results} . \\~\\
Next result is due to Serre :\\~\\
\textbf{Proposition 3.5 :} Assume that the following holds :\\
i) $V_p(\mathfrak{F})$ is a semi-simple $H$ module.\\
ii) $\text{End}(\mathfrak{F}) = \mathbb{Z}_p$.\\
Then, $H_{\text{alg}} = \text{Gl}_h(\mathbb{Q}_p)$ and $H$ is an open sub-group of $\text{Gl}_h(\mathbb{Q}_p)$.\\~\\
\textbf{Proof :} See \cite[Section-5, Theorem-4]{serre 2} . $\square$\\~\\
Now one would like to know when condition-(i) of proposition-3.5 holds. In this direction we have the following result :\\~\\
\textbf{Proposition 3.6 :} The following are equivalent :\\
i) $V_{p}(\mathfrak{F})$ is a semi-simple $H$ module.\\
ii) $V_{p}(\mathfrak{F})$ is a semi-simple $\mathfrak{h}$ module.\\
iii) $\mathfrak{h}$ is a reductive Lie algebra ie a product of an abelian and a semi-simple Lie algebra.\\~\\
\textbf{Proof :} See \cite[Section-1, Proposition-1]{serre 2}. $\square$\\~\\
\textbf{Remark 3.7 :} $\mathfrak{h}$ reductive $\implies$ $\mathfrak{h} = \mathfrak{c} \times \mathfrak{s}$ where $\mathfrak{c}$ is abelian, $\mathfrak{s}$ is semi-simple. In particular, $[\mathfrak{h}, \mathfrak{h}] = \mathfrak{s}$ and $[[\mathfrak{h}, \mathfrak{h}], [\mathfrak{h}, \mathfrak{h}]] = \mathfrak{s}$. \\~\\
To use proposition-3.5 one would like to check if $\mathfrak{h}$ is reductive. We make two definitions and note down some preliminary observations. \\~\\
\textbf{Definition 3.8 :} Let $G$ be a pro-finite group. Put, $G^{(0)} = G$ and $G^{(i)} = [G^{(i-1)}, G^{(i-1)}]$ for each $i \in \mathbb{N}$. \\
i) $G$ is said to be \emph{pro-solvable} if given any open subgroup $U$, $G^{(n)} \subseteq U$ for large enough $n$.\\
ii) $G$ is said to be \emph{almost semi-simple} if $G^{\text{ab}} := G/G^{(1)}$ is finite.\\~\\
\textbf{Remark 3.9 :} i) It is well-known that any finite extension of $E$ is solvable. Hence $\text{Gal}(D_p(\mathfrak{F})|E)$ is pro-solvable. By remark-3.1(iii) $H$ is pro-solvable.\\
ii) Let $G$ be a $p$-adic Lie group over a local field and $\mathfrak{g}$ be the associated Lie algebra. Then $G$ is almost semi-simple if and only if $\mathfrak{g} = [\mathfrak{g}, \mathfrak{g}]$. \\
\subsection{Applications to theory of unramified group laws}
We shall apply the theory described in this section to understand the Galois representation arising from unramified group laws.\\
Let $\mathfrak{F} \in \mathcal{F}_{\pi}(O_E)(h)$. Clearly, $\mathfrak{F} \in \mathcal{F}_p(O_E)(eh)$ where $e$ is the ramification index of the extension $K/\mathbb{Q}_p$. Further, $D_{\pi}(\mathfrak{F}) = D_p(\mathfrak{F})$. Thus, to understand $\text{Gal}(D_{\pi}(\mathfrak{F})|E)$ it is enough to employ techniques described earlier in this section. But to get better results it is desirable to consider representation on $\pi$-adic Tate module $T_{\pi}(\mathfrak{F})$ (to be denoted $\rho_{\pi}(\mathfrak{F})$). For this purpose one needs to improve some of the results used in this section to the context of $p$-adic Lie groups over $K_{\pi}$ and `Hodge-Tate' representations over $K_{\pi}$. A rigorous treatment of this requires lots of effort and so we shall restrict ourselves to the special case $e = 1$ and $\pi = p$. \\
Recall that, the representation $\rho_{\pi}(\mathfrak{F})$ was studied in \cite[Section-4.2]{some results} and the author posed several questions concerning the image $H$ (see \cite[Remark - 4.2.1]{some results}). We make some observations about these questions in the special case $e=1, \pi=p$ using the theory presented in this section. For simplicity, we shall assume that $\mu_{p^h - 1} \subseteq E$. Under this hypothesis, $\text{End}_{O_E}(\mathfrak{F}) = \text{End}(\mathfrak{F})$ (corollary-2.8). \\~\\
\textbf{Remark 3.1.1 :} i) Remark-3.4 shows that $H$ has finite index in an algebraic group (namely $H_{\text{alg}}$) though $H$ may not be algebraic. This gives partial answer to a question in \cite[Remark-4.2.1]{some results}.\\
ii) In the same remark, the author asked a question about minimum dimension of a sub-variety of $\text{Gl}_{h_r}(K_{\pi})$ containing $H$ ($h_r$ as in introduction). Since $H$ is a $p$-adic Lie group the right parameter to study is its dimension as $p$-adic manifold. But by remark-3.4 $H_{\text{alg}}$ has same dimension as $p$-adic manifold. So one can refine the question asked earlier and phrase it in terms of algebraic dimension of $H_{\text{alg}}$.\\
iii) Proposition-3.5 gives answer to another question about $H$ being open, provided the hypothesis on $\mathfrak{h}$ can be verified.\\~\\
Now we make some sketchy remarks about improvement of set-up in general case.\\~\\
\textbf{Remark 3.1.2 :} i) An arbitrary closed subgroup of $\text{Gl}_{h_r}(K_{\pi})$ may not have a $p$-adic manifold structure over $K_{\pi}$. But, since $H \cong \text{Gal}(D_{\pi}(\mathfrak{F})|E)$ and $\mathfrak{F}$ has a $A_{\pi}$ module structure, one expects $H$ to have a $K_{\pi}$-structure. \\
ii) One needs to improve the theory of `Hodge-Tate representations' to vector spaces over $K_{\pi}$ (see \cite[Section-4]{sen} for definition). Note that, one may need to replace cyclotomic character by Lubin-Tate character over $K_{\pi}$. \\
iii) In \cite{sen} Sen uses results of Tate-Sen theory concerning ramifications in $p$-adic Lie extensions. One may like to generalize these results when the Galois group is a $p$-adic Lie group over general local fields. The author intends to revisit this topic in a future article.\\
iv) In \cite{some results} the author has already generalized the concept of $p$-divisible groups to $\pi$-divisible groups in context of $A_{\pi}$ formal modules following \cite{tate} (recall that, all connected $p$-divisible groups arise from divisible formal group laws). One may like to  
prove all the key results of \cite{tate} in this generalized set-up. Note that, if $\mathfrak{F}, \mathfrak{G}$ are connected $p$-divisible groups of dimension 1 both arising from $A_{\pi}$ modules then one can easily generalize the `main result' of $p$-divisible groups to the set-up of $\pi$-divisible groups using the result of Hazewinkel quoted in intoduction.\\
v) We have restricted our attention to 1-dimensional formal groups in the whole theory though many of the results are valid for higher dimensions with slight modifications. The reason behind this restricted approach is inability to get a good ramification theoretic description of torsion subgroups which in 1-dimensional case is provided by Eisenstein polynomials. (see \cite[Appendix-A, Section-2]{the}).
\section{Concluding remarks}    
Let $\mathfrak{F} \in \mathcal{F}_{\pi}(O_E)(h)$ and assume that $h | f_E$. The goal of this section is to pin down some facts about $\text{Gal}(D_{\pi}(\mathfrak{F})|E)$ which will lead to better understanding of the group. In particular, for $e = 1, \pi = p$ case we would like to gather some information towards the problem of finding dimension and properties of $\mathfrak{h}$.\\
For simplicity we shall denote $\text{Gal}(D_{\pi}(\mathfrak{F})|E)$ by $G$ and short-hand $E_{\mathfrak{F}}(\pi^n)$ by $E(\pi^n)$. \\~\\
\textbf{Remark 4.1 :} i) We know that, $\text{Gal}(E(\pi)|E)$ is abelian. Hence $[G, G] \subseteq \text{Gal}(E(\pi^\infty)|E(\pi))$. \\
ii) The $\pi$-adic Tate module $T_{\pi}(\mathfrak{F})$ is a free $A_{\pi}$ of rank $h_r$. Fix a base $z_1, \cdots, z_{h_r}$ of $T_{\pi}(\mathfrak{F})$. Let $z_i(n)$ be the $n$-th component $z_i$ for all $n \in \mathbb{N}$ and $ 1 \leq i \leq h_r$. Put $S_{n} = \{ z_1(n), z_2(n), \cdots z_{h_r}(n) \}$, $\mathfrak{F}'[\pi^n] = \mathfrak{F}[\pi^n] - \mathfrak{F}[\pi^{n-1}]$. Clearly \[K(\mathfrak{F}'[\pi^n]) = K(S_n) = K_{\mathfrak{F}}(\pi^n).\] Let $m_{\pi}(n)$ denote the smallest possible size of a subset of $\mathfrak{F}[\pi^n]$ which generates $K_{\mathfrak{F}}(\pi^n)$ over $K$. Without loss of generality one can assume that such a generating set is a subset of $\mathfrak{F}'[\pi^n]$. Clearly $m_{\pi}(n) \leq h_r$. Consider the sequence $\{m_{\pi}(1), m_{\pi}(2), \cdots \}$. One can conjecture following properties :\\
a) $m_{\pi}(1) \leq  m_{\pi}(2) \leq m_{\pi}(3) \leq \cdots $,\\
b) the sequence is eventually constant and this constant depends only on degree of $F/K_{\pi}$ where $F$ is the fraction field of $\text{End}_{O_E}(\mathfrak{F})$.\\
iii) $m_{\pi}(1) = 1$ and $[F : K_{\pi}] = h_r$ implies $m_{\pi}(i) = 1$ for each $i \geq 1$.\\
Further if $m_{\pi}(i) = 1$ for some $i \geq 2$ looking at degree and ramification index of the extension $K_{\mathfrak{F}}(\pi^i)/K$ we conclude that $K(z) = K_{\mathfrak{F}}(\pi^i)$ for all $z \in \mathfrak{F}'[\pi^i]$. Hence if $m_{\pi}(i) = 1$ for all $i \geq 2$ one has $[F:K_{\pi}] = h_r$ ie $\mathfrak{F}$ has endomorphism ring of full height.  (See \cite[Appendix-A]{the} and \cite[Theorem 2.3]{some results}) \\~\\
The author believes that to get better understanding of $\mathfrak{h}$ one should study questions along these lines and he intends to develope this point of view in a future article.  

\section{Acknowledgement}
I am thankful to Ananyo Kazi for pointing out Sen's work `Ramification in $p$-adic Lie extensions' (1972).

\begin{flushleft}
Kolkata, India\\
(First Version)\\
soumyadip.sahu00@gmail.com
\end{flushleft}


\begin{thebibliography}{wide lebel}

\bibitem[SS1]{some results}
Soumyadip Sahu,
\textit{Field Generated by Division Points of Certain Formal Group Laws;}
2018, arxiv : 1809.00112.

\bibitem[SS2]{the}
Soumyadip Sahu,
\textit{Points of Small Height in Certain Non-Abelian Extensions;}
M.Sc. Thesis, Chennai Mathematical Institute, 2018. arXiv:1806.06587


\bibitem[Haz]{haz2}
M.Hazewinkel;
\textit{Formal Groups and Applications;}
Academic Press, 1978.

\bibitem[Neu]{neu}
J.Neukirch;
\textit{Class Field Theory;}
Springer.


\bibitem[Lub1]{lubin1}
J. Lubin;
\textit{One-Parameter Formal Lie Groups Over p-Adic Integer Rings;}
Annals of Mathematics,  1964.

\bibitem[Lub2]{lubin2}
J. Lubin;
\textit{Finite Subgroups and Isogenies of One Parameter Lie Groups;}
Annals of Mathematics, 1967. 

\bibitem[Ser1]{serre 1}
J.P. Serre;
\textit{Lie Algebras and Lie Groups;}
Springer-Verlag, 2nd Edition.

\bibitem[Ser2]{serre 2}
J.P.Serre;
\textit{Sur les groupes de Galois attach\'es aux groupes p-divisibles;}
Driebergen,1966.

\bibitem[Tate]{tate}
J.T.Tate;
\textit{p-divisible Groups;}
Driebergen, 1966.

\bibitem[Sen]{sen}
S.Sen;
\textit{Lie algebras of Galois groups arising from Hodge-Tate modules;}
Annals of mathematics, 1973.


\bibitem[Wat]{water}
W.C.Waterhouse;
\textit{On p-divisible groups over complete valuation rings}
Annals of mathematics, 1972.

\end{thebibliography}
\end{document}